\documentclass[12pt,titlepage]{article}
\usepackage{amssymb}
\input epsf
\newcommand{\proof}{{\noindent \bf Proof. }}
\newtheorem{thm}{Theorem}
\newtheorem{defi}{Definition}

\newtheorem{lem}{Lemma}[section]

\def\2{\mathbb Z_2}

\title{Relations between the local chromatic number and its directed version}

\author{
      {\bf G\'abor Simonyi}\thanks{Research partially supported by the
Hungarian Foundation for Scientific Research Grant (OTKA) Nos.\ K76088,
K104343, and NK78439, and also by the grant T\'AMOP -
4.2.2.B-10/1-2010-0009.}\\
 Alfr\'ed R\'enyi Institute of Mathematics\\ 
 Hungarian Academy of Sciences\\
 1364 Budapest, POB 127, Hungary\\
 {\tt simonyi.gabor@renyi.mta.hu}
\and
   {\bf G\'abor Tardos}\thanks{Research partially supported by the
NSERC grant 611470 and the
Hungarian Foundation for Scientific Research Grant (OTKA) Nos.\ T037846,
T046234, AT048826, and NK62321.}\\
 Alfr\'ed R\'enyi Institute of Mathematics \\
 Hungarian Academy of Sciences\\
 1364 Budapest, POB 127, Hungary\\
 {\tt tardos.gabor@renyi.mta.hu}
\and
   {\bf Ambrus Zsb\'an}\thanks{Research partially supported by the Hungarian
     Foundation for Scientific Research Grant,  
by the National
Office for Research and Technology (Grant number OTKA 67651). It is also
supported by the grant T\'AMOP - 4.2.2.B-10/1-2010-0009.}\\  
 Department of Computer Science and Information Theory\\
 Budapest University of Technology and Economics \\
{\tt ambrus@math.bme.hu}} 

\date{}

\begin{document}
\maketitle

\begin{abstract}
The local chromatic number is a coloring parameter defined as the minimum
number of colors that should appear in the most colorful closed neighborhood
of a vertex 
under any proper coloring of the graph. Its directed version is the
same when we 
consider only outneighborhoods in a directed graph. For digraphs with all
arcs being present in both directions the two values are obviously equal. Here
we consider oriented graphs. We show the
existence of a graph where the directed local chromatic number of all oriented
versions of the graph is strictly less than the local chromatic number of
the underlying undirected graph. 
We show that for fractional versions the analogous problem has a
different answer: there always exists an orientation for which the directed
and undirected values coincide. We also determine the
supremum of the possible ratios of these fractional parameters, which turns
out to be $e$, the basis of the natural logarithm.
\end{abstract}

\section{Introduction}

The local chromatic number of a graph, first considered in \cite{EFHKRS}, is the
minimum number of colors that must appear in the most colorful closed
neighborhood of a vertex in any proper coloring. Here {\em closed
neighborhood} of a vertex $v$ means $\{v\}\cup N(v)$, where $N(v)$ is the set
of neighbors of $v$. The number of colors is not
restricted, it can be much more than the chromatic number. Formally, denoting
the local chromatic number by $\psi(G)$, we have $$\psi(G):=\min_c \max_{v\in
  V(G)} |\{c(u): u \in N(v)\}|+1,$$ where the minimum is taken over all proper
vertex-colorings $c$ of $G$. Though at first sight one may wonder if this
parameter can ever be smaller than the chromatic number itself, it is shown in
\cite{EFHKRS} that there are graphs with local chromatic number $3$ and
arbitrarily large chromatic number. (It is obvious that $\psi(G)\le \chi(G)$
always holds, where $\chi(G)$ is the chromatic number. It is also easy to see,
that $\psi(G)=2$ implies $\chi(G)=2$.) 

A generalization to directed graphs was defined in \cite{KPS}. The directed
local chromatic number $\psi_d(D)$ of digraph $D$ is defined by 
$$\psi_d(D):=\min_c\max_{v\in V(D)}\{c(u): u\in N_+(v)\}+1,$$
where $c$ runs over all proper vertex colorings of $D$ and $N_+(v)$ is the set
of out-neighbors of $v$, i.e., the set of vertices $w$ with an arc of $D$
going from $v$ to $w$. (A proper vertex-coloring of a digraph is meant to be a
proper vertex coloring of the underlying undirected graph).

It is immediate from the definitions that for a directed graph $\vec G$ and
its underlying undirected graph $G$ we always have $\psi_d(\vec G)\le \psi(G)$
and we automatically have equality if $\vec G$ contains both orientations of
all edges of $G$.

We say that $\vec G$ is an orientation of $G$ if $\vec G$ contains all
edges of $G$ with exactly one of the two possible orientations and no other
arcs. Whether we can
achieve the value $\psi(G)$ as the directed local chromatic number of an
orientation of $G$ is a natural question, and was already
asked in \cite{STjgt}. We will show in this paper that this is not always the
case.

With regard to the fractional versions of these chromatic parameters
(also introduced in \cite{KPS}), we show that every undirected graph $G$ has
an orientation the fractional directed local chromatic number of which
achieves the fractional local
chromatic number of $G$. The latter happens to coincide with the fractional
chromatic number as shown in \cite{KPS}. We also find the supremum of the
fractional (local) chromatic number of the underlying graph of any
digraph with a given fractional directed local chromatic number.

Observe the analogy with Shannon capacity of graphs and its generalization,
the Sperner capacity of digraphs, cf. \cite{GKV3}. (We mention the connection
that the Sperner capacity is always bounded from above by the fractional
directed local chromatic number, see \cite{KPS}.) The Sperner capacity of a
digraph $\vec G$ is at most the Shannon capacity of the underlying graph $G$
and equality is achieved if both orientations of every edge of $G$ is present
in $\vec G$.
Whether the Shannon capacity of an undirected graph can always be achieved as
the Sperner capacity of an orientation
is an open question. It was investigated in \cite{SaSi}, where it was
shown that this is always the case for a non-trivial class of graphs. In the
light of that result and no indication that this is
not always true, the negative answer in case of the local
chromatic number is surprising.

The paper is organized as follows. In Section \ref{gapegy} we show that
there exists a graph $G$ such that all orientations $\vec G$ of $G$ satisfy
$\psi_d(\overrightarrow{G})<\psi(G)$. 
A lemma (a special case of which is) needed to proving this statement has a
somewhat tedious proof, therefore the proof of this lemma is postponed to
Section \ref{lemmak}.
Section \ref{fraceq} is devoted to the analogous question for
fractional versions of the above parameters. We show that, in contrast to the
result in Section \ref{gapegy}, the analogous
inequality is always an equality in the fractional case. In Section
\ref{bigrat} we 
give tight bounds for the largest possible ratio of these fractional
parameters. As already mentioned, Section \ref{lemmak} contains the proof of 
Lemma \ref{unicolor}.

We note that the question answered in Section \ref{bigrat} recently turned out
to be relevant also for a problem of information transmission, see \cite{SDL},
where a somewhat weaker upper bound was proved independently.

\section{Local chromatic number: the oriented and undirected case}\label{gapegy} 

For a graph $G$ let
$$\psi_{d,\max}=\max\{\psi_d(\vec{G}):
  \vec{G} {\rm \ is\ an\ orientation\ of\ } G\}.$$
Here we will construct a graph proving the following statement.

\begin{thm} \label{gap1}
There exists a graph $G$ such that $$\psi_{d,\max}<\psi(G).$$ 
\end{thm}

To prove this theorem, we will consider the graphs $U(m,k)$ that are the
universal graphs for the local chromatic number in the following sense: 
A graph $G$ has a local $k$-coloring with at most $m$ colors if and only if
$G$ admits a homomorphism to $U(m,k)$. Here a {\em local $k$-coloring} is a
proper coloring with the vertices of any closed neighborhood receiving at most
$k$ distinct colors.

\begin{defi} \label{Unk} {\rm (\cite{EFHKRS})}
Let $k\leq m$ be positive integers and $[m]=\{1,2,\ldots,m\}$. The graph
$U(m,k)$ is defined as follows.
$$V(U(m,k)):=\{(x,A): x\in [m], A\subseteq [m], |A|=k-1, x\notin A\}$$ and
$$E(U(m,k)):=\{\{(x,A),(y,B)\}: x\in B, y\in A\}.$$
\end{defi}

The {\em natural coloring} of
$U(m,k)$ colors the vertex
$(x,A)\in V(U(m,k))$ with color $x$. It is easy to see that this is a local
$k$-coloring.

For the proof of Theorem \ref{gap1} we will use the fact that the natural
coloring is basically the only local
3-colorings of $U(5,3)$. We will prove this statement 
more generally. 
 
\begin{lem}\label{unicolor}
Let us have $m>k+1\ge 4$ and let $c$ be a local
$k$-coloring of the graph $U(m,k)$. Then $c$ is the natural coloring up to
permutation of colors.
\end{lem}

As the proof of this lemma is a little tedious, we postpone it
to Section \ref{lemmak}. Here we present the proof of Theorem \ref{gap1}
assuming that Lemma \ref{unicolor} is true for m=5 and k=3.

\medskip
\par\noindent
{\bf Proof of Theorem \ref{gap1}}
We construct a graph $G$ that has the property given in the statement.
Let $$V(G)=V(U(5,3))\cup \{x,y,z\},$$
with $x,y,z\notin V(U(5,3))$, while the edge set
is $$E(G)=E(U(5,3))\cup\{\{x,y\},\{x,z\},\{y,z\}\}\cup\{\{x,(2,\{1,3\})\}, 
\{y,(3,\{1,2\})\}, \{z,(1,\{4,5\})\}.$$ 
\par\noindent 
First we show that $\psi(G)\ge 4$ (in fact it is exactly $4$). Assume for
contradiction that $f$ is a local 3-coloring of $G$. Notice that $U(5,3)$ is a
subgraph of $G$, so by Lemma~\ref{unicolor}
we may assume that $f$ is the natural coloring on $U(5,3)$, i.e.,
we assume $f((i,A))=i$ for every choice of $(i,A)\in V(U(5,3))$. 
\par\noindent 
The above means that $f(x)$ should be $1$ or $3$ as otherwise we have too many
colors in the neighborhood of $(2,\{1,3\})$. Similarly $f(y)$ should be $1$
or $2$ and $f(z)$ should be $4$ or $5$.
\par\noindent 
If $f(x)=1$ then $y$ is connected to vertices of color $1$, $3$ and $4$ or $5$,
too many distinct colors.
\par\noindent
If $f(y)=1$ then $x$ is connected to vertices of color $1$, $2$ and $4$ or $5$,
which is also too many.
\par\noindent
The only remaining possibility is $f(x)=3$ and $f(y)=2$, but then the
neighborhood of $z$ has too many colors.
\par\noindent
The contradiction proves $\psi(G)\ge 4$. 
\medskip
\par\noindent
Now we show that $G$ has no orientation $\vec{G}$ that has directed local
chromatic number more than $3$.
\medskip
\par\noindent
We will focus only on the orientation of the edge $\{x,y\}$. If it is oriented
from $x$ to $y$, then the following coloring $g$ of $\vec{G}$ will show
$\psi(G)\le3$. Let $g(x)=1$, $g(y)=2$ while $U(5,3)$ gets its natural coloring
and $g(z)$ is either $4$ or $5$. Then the closed outneighborhood of no vertex
contains more than $3$ colors. For the opposite case when the $\{x,y\}$ edge
is oriented from $y$ to $x$ we define the coloring $g'$. Let $g'(x)=3,
g'(y)=1$ and $g'(v)=g(v)$ for all vertices $v\neq x,y$. 
It is easy to check that no closed outneighborhood contains more
than $3$ colors in this case either. So 
$\psi_{d,\max}(G)\le3$ is strictly smaller than $\psi(G)\ge4$ and the
proof is complete.
\hfill$\Box$
\medskip
\par\noindent
It is annoying that we do not know an example where the gap between $\psi(G)$
and $\psi_{d,\max}(G)$ is more than $1$.

\section{The maximum fractional directed local chromatic number of an orientation}\label{fraceq}

Unlike the case of the integral values discussed in the previous section,
there always exists an orientation of any graph for which the 
fractional 
relaxation of the directed local chromatic number attains
the fractional (undirected) local chromatic number (which is just the
fractional chromatic 
number, see \cite{KPS}). This is what we prove in this section. 

\subsection{Linear programming definitions}

All fractional graph parameters are defined as the optimum value of certain
linear programs. We start by recalling the definition of fractional
colorings and the fractional chromatic number.

Let $G$ be a graph. Let $S(G)$ denote the set of
independent sets of $G$. A {\em fractional coloring} of $G$ is the collection
of real weights $x_A$ for independent sets $A\in S(G)$ satisfying
the following linear inequalities. 
$$\forall A:\ x_A\ge 0$$
\begin{equation}
\forall v\in V(G):\ \sum_{v\in A\in S(G)} x_A \ge 1
\end{equation}

The fractional chromatic number $\chi^*(G)$ of $G$ is the minimal total weight
of a fractional coloring, that is:
$$\chi^*(G)=\min\sum_{A\in S(G)}x_A,$$
where the minimum is taken for all fractional coloring $(x_A)$ of $G$. (We can
write minimum here as it is attained in this and similar LP problems.)

Note that a proper coloring of $G$ can be turned into a fractional coloring
by giving weight $1$ to the color classes and weight $0$ to all other
independent sets. This shows $\chi^*(G)\le\chi(G)$, where the chromatic number
$\chi(G)$ of $G$ is the smallest number of colors in a proper vertex coloring
of $G$.

Let $\vec G$ be a directed graph and $G$ the underlying undirected graph. By a
fractional coloring of $\vec G$ we mean a fractional coloring of $G$. The
{\em local weight} of a fractional coloring $(x_A)$ is
$$1+\max_{v\in V(\vec G)}\sum_{A\in S(G):N_+(v)\cap A\ne\emptyset}x_A,$$
i.e., 1 plus the maximum total weight a vertex ``sees''. The
fractional directed local chromatic number $\psi_d^*({\vec G})$ is
then the minimal local weight of a fractional coloring of $\vec G$.

Note that if we turn a directed local $k$-coloring of $\vec G$ into a
fractional coloring it has local weight at most $k$, thus we have
$\psi_d^*(\vec G)\le\psi_d(\vec G)$ as expected. Here a {\em directed local
$k$-coloring} is a proper coloring with every out-neighborhood receiving at
most $k-1$ colors.

The {\em fractional local chromatic number} $\psi^*(G)$ of a graph can be
defined as $\psi_d^*(\vec G)$, where $\vec G$ is obtained from $G$ by
replacing each of its edges by the two arcs representing its two
orientations. Note that this is not new graph parameter, but rather we have
$\psi^*(G)=\chi^*(G)$ as proved in \cite{KPS}.

To compare the fractional directed local chromatic number of a digraph to the
fractional chromatic number of the underlying undirected graph, we will use
dual formulation of latter as the fractional clique number. That is, we use
that the fractional chromatic number of a graph is the maximum total weight
assigned to the vertices satisfying that (i) all weights are non-negative and
(ii) the total weight of the vertices of an independent set does not exceed
$1$.

\subsection{Equality for fractional values}

We have already mentioned the result $\psi^*(G)=\chi^*(G)$ from
\cite{KPS}. This statement claims that with any fractional coloring of a graph
$G$ there will be a vertex $v\in V(G)$ that ``sees'' a total weight at least
$\chi^*(G)-1$. Here we prove that $v$ can be chosen from a fixed independent
set.

\begin{thm}\label{ize}
For a graph $G$ and a vertex $v_0\in V(G)$ there exists an independent set
$A_0\in S(G)$ containing $v_0$ such that the following holds. For any
fractional coloring $(x_A)_{A\in S(G)}$ of $G$ there is a vertex $v\in A_0$ such
that
$$\sum_{A\in S(G): N(v)\cap A\ne\emptyset}x_A\ge\chi^*(G)-1.$$ 
\end{thm}

\proof
Consider an optimal fractional clique $(t_v)_{v\in V(G)}$ of $G$, that
is one attaining
\begin{equation}
\sum_{v\in V(G)}t_v=\chi^*(G).
\end{equation}
Note that $t_v\ge0$ for all $v\in V(G)$ and
\begin{equation}
\sum_{A\ni v}t_v\le1
\end{equation}
for every independent set $A$ of $G$. Let us choose an independent
set $A_0$ containing $v_0$ for which
\begin{equation}
\sum_{i\in A_0} t_i = 1.
\end{equation}
Such a set $A_0$
must exist, otherwise the value of $t_{v_0}$ could be increased showing
that the fractional clique $(t_v)_{v\in V(G)}$ does not have maximal total
weight.

Now let $(x_A)_{A\in S(G)}$ be any fractional coloring of $G$. We have
$$\begin{array}{rll}
\chi^*(G)-1=&\sum_{v\notin A_0}t_v&\hbox{by (2) and (4)}\\
\le&\sum_{v\notin A_0}t_v\sum_{A\ni v}x_A&\hbox{by (1)}\\
=&\sum_{A\in S(G)}x_A\sum_{v\in A\setminus A_0}t_v\\
\le&\sum_{A\in S(G)}x_A(1-\sum_{v\in A_0:N(v)\cap A=\emptyset}t_v)&\hbox{see below}\\
=&\sum_{A\in S(G)}x_A\sum_{v\in A_0:N(v)\cap A\ne\emptyset}t_v&\hbox{by (4)}\\
=&\sum_{v\in A_0}t_v\sum_{A:N(v)\cap A\ne\emptyset}x_A\\
\le&\max_{v\in A_0}\sum_{A:N(v)\cap A\ne\emptyset}x_A&\hbox{by (3),}
\end{array}$$
where the marked inequality follows from (3) applied to the independent set
$A\cup \{v\in A_0\mid N(v)\cap A=\emptyset\}$.

Comparing the first and last lines above proves the theorem.
\hfill$\Box$ 

\begin{thm}\label{frakceq}
Let $G$ be a finite undirected
graph. Then $$\max_{\vec G}\psi_d^*(\vec G)=\chi^*(G),$$ 
where the maximum is taken over all orientations $\vec G$ of $G$.
\end{thm}

\proof
It is clear that $\chi^*(G)=\psi^*(G)$ is an upper bound on the left hand
side. We
need to give an orientation $\vec G$ of $G$ satisfying $\psi_d^*(\vec
G)\ge\chi^*(G)$.

Consider any independent set $A_0$ of $G$ that satisfies the
statement of Theorem~\ref{ize}. Let us obtain $\vec G$ from $G$ by orienting
each edge connecting $A_0$ to its complement in the direction leaving $A_0$
and orienting the remaining edges of $G$ arbitrarily.

Now consider any fractional coloring $(x_A)_{A\in S(G)}$ of $\vec G$. This is
a fractional coloring of $G$ and by the choice of $A_0$ there is a vertex
$v\in A_0$ with $\sum_{A:N(v)\cap A\ne\emptyset}x_A\ge\chi^*(G)-1$. Here $N(v)$
refers to neighborhood of $v$ in $G$, but the orientation we chose makes this
equal to the outneighborhood $N_+(v)$ of $v$ in $\vec G$. Thus we also have
$\sum_{A:N_+(v)\cap A\ne\emptyset}x_A\ge\chi^*(G)-1$ finishing the proof.
\hfill$\Box$ 

\section{The smallest fractional directed local chromatic number}\label{bigrat}

In this section we determine the supremum of $\chi^*(G)$ given a fixed bound
on $\psi_d^*(G)$ and, as a consequence, we also find the supremum of
the ratio 
$\frac{\chi^*(G)}{\psi^*_d(G)}$.  Here $G$ is a directed graph (note that we
dropped the arrow from the notation $\vec G$ used in the earlier
sections). Recall that $\chi^*(G)$ and $\psi^*_d(G)$ are 
the fractional chromatic number and the fractional directed local chromatic
number, respectively, and the former is defined as the fractional
chromatic number of the underlying undirected graph, and is also equal to the
fractional local chromatic number of this undirected graph (\cite{KPS}). We
also note that the boundedness of the above ratio was independently proved in
\cite{SDL}, where the somewhat weaker upper bound, $\frac{5}{4}e^2$, was
presented.

\begin{thm}\label{ratio} 
\par\noindent
(a) For every finite, loopless directed graph $G$ we have
$$\chi^*(G)\le\frac{k^k}{(k-1)^{k-1}}<ek,$$
where $k=\psi^*_d(G)>1$ and $e$ is the basis of the natural
logarithm. 

(b) For every $k\ge2$ and $\varepsilon>0$ there exists a loopless finite
directed 
graph $G$ with $\psi^*_d(G)\le k$ and
$$\chi^*(G)>\frac{k^k}{(k-1)^{k-1}}-\varepsilon.$$
If $k$ is an integer, then the above graph can be chosen to further satisfy
$\psi_d(G)=k$. 
\end{thm}

Note that for graphs $G$ with no edges we have $\chi^*(G)=\psi_d^*(G)=1$,
while non-trivial graphs have $\psi_d^*(G)\ge2$. Therefore Theorem~\ref{ratio}
establishes these tight results:
\begin{eqnarray*}
\sup\{\chi^*(G)\mid\psi^*_d(G)\le
  k\}&=&\left\{\begin{array}{cl}\frac{k^k}{(k-1)^{k-1}}&\mbox{ for
      }k\ge2\\1&\mbox{ for }1\le k<2,\end{array}\right.\\
\sup\frac{\chi^*(G)}{\psi^*_d(G)}&=&e.
\end{eqnarray*}

Before proving the theorem above we give simple alternative definitions for
both graph parameters concerned.

\begin{lem}\label{fracchrom}
Let $G=(V,E)$ be a (directed or undirected) graph. We have
$$\chi^*(G)=(\sup\min_{v\in V}P[v\in I])^{-1},$$
where the supremum is over random variables $I$, whose values are independent
sets of $G$.
\end{lem}

\proof
The fractional chromatic number $\chi^*(G)$ is defined as the minimum of the
total weight $s$ of fractional colorings $(x_A)_{A\in S(G)}$ of $G$. The $\le$
direction of the lemma is proved by considering
the random variable $I$ that takes an independent set $A$ with probability
$x_A/s$. Now let $I$ be a random variable taking values
from $S(G)$ and let $c=\min_{v\in V}P[v\in I]$. For the reverse direction
consider the fractional coloring given by $x_A=P[I=A]/c$.
\hfill$\Box$

\medskip
\par\noindent
Let $G=(V,E)$ be a directed graph, $C$ be an arbitrary finite set (the set of
colors) and $r\ge1$ an integer. We denote by $C\choose r$ the family of
subsets of $C$ of size $r$. We call a function $\chi:V\to {C\choose r}$ an
{\em$r$-multi-coloring} of $G$ if for all $c\in C$ the set $\{v\in V\mid c\in
f(v)\}$ is independent in $G$. Note that 1-multi-coloring is a proper
coloring and in general $r$-multi-coloring is a homomorphism to the appropriate
Kneser graph. An $r$-multi-coloring $\chi$ of $G$ is an {\em$h$-local
$r$-multi-coloring} if $|\bigcup_{w\in N^+(v)}\chi(w)|\le h-r$ for all $v\in V$.
The following easy lemma was already used in \cite{KPS}.

\medskip\par\noindent
\begin{lem}\label{multi}
The fractional directed local chromatic number $\psi_d^*(G)$ of a directed
graph $G$ is the infimum of the fractions $h/r$ such that a $h$-local
$r$-multi-coloring of $G$ exists.
\end{lem}

\medskip\par\noindent
{\bf Proof of Theorem~\ref{ratio}.}
(a)
In light of Lemmas~\ref{fracchrom} and \ref{multi} to prove the upper bound on
the fractional chromatic number it is enough consider an $h$-local
$r$-multi-coloring $\chi$ of a directed graph $G$ and define a random variable
$I$ taking independent sets as values that satisfies
$$P[v\in I]\ge\frac{(h/r-1)^{h/r-1}}{(h/r)^{h/r}}$$
for every vertex $v\in V$.

Let us ``select'' each of the colors used by $\chi$ independently and with the
same probability $1-\gamma$ to be set later. Let $C'$ stand for the set of
these selected colors. Let $I=\{v\in V\mid \chi(v)\cap C'\ne\emptyset,\forall
w\in N^+(v):\chi(w)\cap C'=\emptyset\}$. Here $I$ consists of the vertices
$v\in V$ with at least one of their $r$ colors selected but satisfying that
none of their out-neighbors have any selected colors. Clearly, this
is an independent set. Let us fix a vertex $v$. We have $v\in I$ if one of the
colors in $\chi(v)$ is selected but none of the colors in $S=\cup_{w\in
N^+(v)}\chi(w)$ is selected. Here $S$ and $\chi(v)$ are disjoint by the
definition of multi-coloring, so these events are independent and we have
$P[v\in I]=(1-\gamma^r)\gamma^{|S|}$. We have $|S|\le h-r$ as $\chi$ is an
$h$-local $r$-multi-coloring, so $P[v\in I]\ge(1-\gamma^r)\gamma^{h-r}$. Setting
$\gamma=(1-r/h)^{1/r}$ gives the desired bound and finishes the proof of the
upper bound.

\medskip
\par\noindent
(b)
Let now $m\ge k$ be positive integers and $C$ a set of $m$ colors. We consider
the directed analogue of the universal graph $U(m,k)$ used above. The directed
graph $U_d(m,k)$ is defined as follows. The vertex set consists of pairs 
$(c,H)$ consisting of $c\in C$ and $H\subset C$ satisfying $c\notin H$ and
$|H|=k-1$. We set a directed edge from $(c,H)$ to $(c',H')$ if $c'\in H$. It
is easy to see that $U_d(m,k)$ is the universal graph for directed local
$k$-coloring using $m$ colors, that is, a directed graph has such a coloring
if and only if it has a homomorphism to $U_d(m,k)$. We will not use this
universality but we use that the {\em natural coloring} $\chi_0$ mapping the
vertex $(x,H)$ to $x$ is a directed local $k$-coloring showing $\psi_d(G)\le
k$. In 
fact, we have equality here since $U_d(m,k)$ contains a complete directed
subgraph on the $k$ vertices $(x,H\setminus\{x\})$ for some $k$-element subset
$H$ of $C$.

To bound the fractional chromatic number of $U_d(m,k)$ we use the simple
observation that for any graph $G=(V,E)$ we have $\chi^*(G)\ge|V|/\alpha(G)$,
where $\alpha(G)$ is the size of the largest independent set in $G$. We need a
lower bound only, but remark that for vertex-transitive graphs such as
$U_d(m,k)$ equality holds.

Let $I$ be a largest independent set of $U_d(m,k)$ and let $C'=\{\chi_0(v)\mid
v\in I\}$ be the set of their colors under the natural coloring. For a vertex
$v=(c,H)\in I$ we must have $c\in C'$ and $H\cap C'=\emptyset$. All the
vertices of $U_d(m,k)$ satisfying this form an independent set of size
$(m-l){l\choose k-1}$, where $l=m-|C'|$. As $U_d(m,k)$ has $m{m-1\choose
  {k-1}}$ 
vertices we have
$$\chi^*(U_d(m,k))\ge\frac{m{m-1\choose k-1}}{\max_l((m-l){l\choose
    k-1})}\ge\frac{m(m-1)^{k-1}}{\max_l((m-l)l^{k-1})}\ge(1-\frac1m)^{k-1}\frac{k^k}{(k-1)^{k-1}}.$$ 

The last inequality above follows from realizing that $l=m-m/k$ maximizes the
denominator. The lower bound obtained here is exactly the desired bound except
the multiplicative error term $(1-1/m)^{k-1}$. The effect of this error term
can be made arbitrarily small by choosing $m$ large enough. This proves the
last statement of the theorem.

It remains to prove the first statement of part (b), namely the tightness of our upper
bound on the fractional chromatic number for graphs with a non-integer
fractional directed local chromatic number. We do this similarly to our proof
for the integer case but have to consider universal graphs for $h$-local
$r$-multi-colorings.

Let $m\ge h$ and $r\le h/2$ be positive integers, let $C$
be a set of $m$ colors and consider the directed graph $U_d(m,h,r)$ whose
vertices are pairs $(Q,H)$ satisfying $Q\in{C\choose r}$, $H\in{C\choose
h-r}$ and $Q\cap H=\emptyset$. We have a directed edge from $(Q,H)$ to
$(Q',H')$ if $Q'\subseteq H$. Clearly, $U_d(m,h,r)$ is the universal graph for
having $h$-local $r$-multi-coloring using $m$ colors. (Note that the
undirected version $U(m,h,r)$ of $U_d(m,h,r)$ is just the graph denoted
$U_r(m,h)$ in \cite{KPS}.) The {\em natural
multi-coloring} maps a vertex $(Q,H)$ to $Q$ and this is clearly an
$h$-local $r$-multi-coloring and shows $\psi_d^*(U_d(m,h,r))\le h/r$. 

Let $I$ be a largest independent set in $U_d(m,h,r)$ and set $\mathcal
H=\{H|\exists Q:(Q,H)\in I\}\subseteq{C\choose h-r}$. Let $\mathcal S$ be the
size $r$ shadow of $\mathcal H$, that is $\mathcal S=\cup_{H\in\mathcal
H}{H\choose r}$. Let us find the value $l$ (not necessarily integer) such that
$|\mathcal H|={l\choose h-r}$. By the general form of the Kruskal-Katona
theorem \cite{LLkonyv} (exercise 13.31)
we have $|\mathcal C|\ge{l\choose r}$. A vertex
$(Q,H)\in I$ 
satisfies $Q\notin\mathcal S$ and $H\in\mathcal H$. Thus we have
$$\alpha(U_d(m,h,r))=|I|\le\left({m\choose r}-{l\choose r}\right){l\choose
  h-r}.$$ 
We use the inequalities ${m\choose r}-{l\choose r}\le(m^r-l^r)/r!$ and
${l\choose h-r}\le l^{h-r}/(h-r)!$ plus calculus to obtain:
$$\alpha(U_d(m,h,r))\le\frac{(m^r-l^r)l^{h-r}}{r!(h-r)!}\le\frac{(h/r-1)^{h/r-1}}{(h/r)^{h/r}}\cdot\frac{m^r}{r!(h-r)!}.$$ 

Finally using that $U_d(m,h,r)$ has $n={m\choose r}{m-r\choose
h-r}\ge(m-h)^h/(r!(h-r)!)$ vertices we obtain
$$\chi^*(U_d(m,h,r))\ge\frac n{\alpha(U_d(m,h,r))}\ge(1-h/m)^h\cdot\frac{(h/r)^{h/r}}{(h/r-1)^{h/r-1}}.$$

To finish the proof of the theorem let $s\ge2$ be an
arbitrary real. Take positive integers $m\ge h$ and $r$ such that $2\le h/r\le
s$. Consider the graph $U_d(m,h,r)$. It satisfies $$\psi_d^*(U(m,h,r))\le h/r\le
s,$$
$$\chi^*(U_d(m,h,r))\ge(1-h/m)^h(h/r)^{h/r}/(h/r-1)^{h/r-1}.$$ If we
choose $h/r$ close enough to $s$ and $m$ large enough this last value will be
arbitrarily close to $s^s/(s-1)^{s-1}$.
\hfill$\Box$

\medskip\par\noindent
{\bf Remark 1.} Part (a) of Theorem \ref{ratio} can also be proven using the
proof method of Theorem 5 in \cite{KPS} stating
$\psi^*(G)=\chi^*(G)$. On the one hand we can use that $G$ having an $h$-local
$r$-multi-coloring with $m$ 
colors is equivalent to $G$ admitting a homomorphism to $U_d(m,h,r)$. On the
other hand we can also use the fact that if $G$ admits a homomorphism to
another graph $H$, then $\chi^*(G)\le \chi^*(H)$. Thus it is enough to
determine the fractional chromatic number of $U_d(m,h,r)$ (and maximize its
value in $m$ while $h/r$ is fixed), to get the largest possible fractional
chromatic number of a graph with fractional directed local chromatic number
$h/r$. Using the vertex-transitivity of the graphs $U_d(m,h,r)$ this can be
done by determining their independence number and using that for a
vertex-transitive graph $F$ the fractional chromatic number $\chi^*(F)$ is
equal to $\frac{|V(F)|}{\alpha(F)}$. This way a more precise calculation in
the proof of part (b) can actually lead to a proof for part (a), too. 
\hfill$\diamond$

\medskip\par\noindent
{\bf Remark 2.} 
As we have seen the fractional chromatic and the fractional local chromatic
numbers agree for undirected graphs, but the former may be larger by a
factor up to $e$ if we consider directed graphs. However, if we consider
another variant of these chromatic parameters for directed graphs,
originally introduced by Neumann-Lara \cite{NL}, then the two parameters agree
for directed graphs as well. For the definitions of these variants one uses
acyclic subsets (i.e., subsets of the vertices that induce acyclic graphs)
the same way as one uses independent sets for the variants we studied.
\hfill$\diamond$

\section{Proof of Lemma \ref{unicolor}}\label{lemmak}

\par\noindent
In this section we do not distinguish colorings that differ only in the
permutations of colors, i.e., we identify colorings that induce the same
partition on the vertex set.

For the proof of Lemma \ref{unicolor} we will need the following other lemma. 
Notice that the graph $U(k+1,k)$ is a $k$-chromatic graph, thus it is not an
example for having strictly larger chromatic number than local chromatic
number. On the other hand, all the graphs $U(m,k)$ with $m>k+1>3$ give such
examples.   

As opposed to $U(m,k)$ with $m>k+1>3$, the graph $U(k+1,k)$ have local
$k$-colorings different from the natural coloring. Indeed any proper
$k$-coloring is also a local $k$-coloring. Our next lemma states that for
$k\ge 3$ there are no other local $k$-colorings. 

\begin{lem}\label{k+1,k}
For every $k\ge 3$ the graph $U(k+1,k)$ has $k+2$ different local
$k$-colorings: its natural coloring and the $k$-colorings obtained from the
natural coloring by keeping the color of all but a single color class and
recoloring the vertices of this last color class the only way to obtain a
proper $k$-coloring.
\end{lem}

\proof Let $f$ be a local $k$-coloring of $U(k+1,k)$.
For simplicity we write $(a,\overline b)$ for the vertex $(a,[k+1]\setminus\{a,b\})$
and $f(a,\overline b)$ for $f((a,\overline b))$. Let $m$ stand for the number of colors $f$
uses.

The vertices $(i,\overline{k+1})$ for $i\in[k]$ form a $k$-clique. As $f$ is proper
it assigns distinct colors to these vertices. We may therefore assume
$f(i,\overline{k+1})=i$ for $i\in[k]$. As $f$ is a local $k$-coloring all the
neighbors of the vertices in this clique must also receive colors from
$[k]$. Thus the $m-k$ colors outside $[k]$ that $f$ uses must appear on
vertices receiving color $k+1$ in the natural coloring (all other vertices
are neighbors of our clique).

By symmetry all the $k+1$ color classes of the natural coloring must
contain $m-k$ color classes of the coloring $f$, so we must have
$m\ge(k+1)(m-k)$. This means that $m=k$ or $m=k+1$. Furthermore, in case
$m=k+1$ all $k+1$ color classes must live inside a single color class of
the natural coloring, meaning that $f$ is equivalent to the natural
coloring.

It remains to consider the case $m=k$, thus $f$ is a proper $k$-coloring.

For $i\ne j$, $i,j\in[k]$ we have $f(i,\overline j)\in\{i,j\}$ as $(i,\overline j)$ is
connected to $(l,\overline{k+1})$ of color $l$ for $l\in[k]\setminus\{i,j\}$. We
call a vertex $(i,\overline j)$ {\em special} if $f(i,\overline j)=j$ (even if $i=k+1$). If a
vertex $(i,\overline j)$ for $i\in[k]$ is not special we have $f(i,\overline j)=i$.

Note that the vertices $(i,\overline j)$ for a fixed $j$ form a clique, thus at most
one of them can be special (having color $j$) or $f$ is not proper.

If all the vertices $(k+1,\overline i)$ are special, then these are all the special
vertices and $f$ is obtained from the natural coloring by recoloring its
last color class to the remaining colors.

Otherwise we have $f(k+1,\overline i)=j$ for some $i,j\in[k]$, $i\ne j$. To make $f$
proper all the vertices $(j,\overline l)$ with $l\in[k]\setminus\{j\}$ must be
special. This means that no more vertices can be special except possibly a
single vertex $(x,\overline j)$ for some $x$. Indeed $(k+1,\overline j)$ must be special,
moreover $f(k+1,\overline l)=j$ for all $l\in[k]$ as $(k+1,\overline l)$  is connected to
vertices of all colors but $j$, namely to the special vertices $(j,\overline{l'})$ if
$l\ne j$ and to the non-special vertices $(l',\overline{l''})$ with
$l',l''\in[k]\setminus\{j\}$ if $l=j$. This makes $f$ equivalent to the
$k$-coloring obtained from the natural coloring by recoloring color class
$j$.
\hfill$\Box$

\medskip
\par\noindent
{\bf Proof of Lemma \ref{unicolor}.} We use the phrase {\em class} for
the color classes of the natural coloring in $U(m,k)$. We write $c(x,A)$ for
the color $c((x,A))$ of the vertex $(x,A)$.

First we consider the $m=k+2$ special case.

For $i\in[k+2]$ consider the
subgraph $G_i$ of $U(m,k)$ induced by the vertices $(x,A)$ with $i\notin
A\cup\{x\}$. Clearly, $G_i$ is isomorphic to $U(k+1,k)$, and thus the
restriction to $G_i$ of the local $k$-coloring $c$ must be one of the few
colorings described in Lemma~\ref{k+1,k}. In particular, the intersection of
a class with $G_i$ either receives a single color by $c$ (a monochromatic
intersection) or each vertex of the intersection receives different colors (a
colorful intersection). Furthermore, $k+1$ classes intersect $G_i$ but at
most one of these intersections is colorful.

Let us fix a color class of the natural coloring of $U(m,k)$. It intersects
$k+1$ of the $k+2$ subgraphs $G_i$. If at least two of these intersections are
monochromatic, then all of them have to be, so $c$ assigns the same color to
the entire class. Otherwise at least $k$ of the intersections are
colorful. Since altogether there are at most $k+2<2k$ colorful
intersections, we must have at most a single class to which $c$ assigns
multiple colors.

If all the classes are monochromatic, then $c$ must assign distinct colors to
them or it is not a proper coloring. This makes $c$ equivalent to the natural
coloring as needed. Otherwise we have a single non-monochromatic class with
$k$ or $k+1$ of its non-empty intersections with the subgraphs $G_i$
colorful. All that other classes are monochromatic and must receive distinct
colors, thus we assume that $c$ equals to the natural coloring outside the
only exceptional class.

Each vertex in this exceptional class is contained in exactly two of the
subgraphs $G_i$, and since more than half these intersections are colorful
there must be a vertex $(x,A)$ that is in two distinct subgraphs $G_i$ and
$G_j$ with the intersection of either of them with the exceptional class being
colorful. Lemma~\ref{k+1,k} determines the color of $(x,A)$ from the coloring
of $G_i$ outside the exceptional class: it must be $j$. The contradiction
finishing the proof of the $m=k+2$ case of the lemma comes from observing that
a similar argument using the restriction of $c$ to $G_j$ gives $c(x,A)=i$.

Finally, let us consider the case $m>k+2$. For any subset $H\subset[m]$ with
$|H|=k+2$ the vertices $(x,A)$ satisfying $A\cup\{x\}\subset H$ induce a
subgraph $G_H$ of $G$ isomorphic to $U(k+2,k)$. Thus $c$ must be equivalent to
the natural coloring on $G_H$. This means that $c(x,A)=c(x,B)$ whenever $|A\cup
B|\le k+1$ as in this case $(x,A)$ and $(x,B)$ are in a common subgraph
$G_H$. Now consider two arbitrary vertices $(x,A)$ and $(x,B)$ from the same
class. Clearly there exists a sequence $A=A_0,A_1,\ldots,A_k=B$ such that
$c(x,A_{i-1})=c(x,A_i)$ for the above reason for all $i\in[k]$. Thus we have
all $c(x,A)=c(x,B)$ and all the classes are monochromatic at $c$. To
make $c$ proper it must assign distinct colors to to distinct classes, thus
$c$ is equivalent to the natural coloring. This finishes the proof of the
lemma.
\hfill$\Box$

\end{document}